\documentclass[10pt,twoside]{amsart}
\usepackage{amsmath}
\usepackage{amssymb}
\usepackage{amscd}
\usepackage{amsthm}
\theoremstyle{plain}
\newtheorem{Thm}{Theorem}[section]
\newtheorem{Prob}[Thm]{Problem}
\newtheorem{Prop}[Thm]{Proposition}
\theoremstyle{definition}
\newtheorem{Defn}[Thm]{Definition}
\theoremstyle{remark}

\numberwithin{equation}{section}
\newcommand{\Exc}{\mathop{\mathrm{Exc}}}
\newcommand{\Supp}{\mathop{\mathrm{Supp}}}
\newcommand{\Strata}{\mathop{\mathrm{Strata}}}
\newcommand{\CNKLT}{\mathop{\mathrm{C}_\mathrm{non-klt}}}
\newcommand{\MCNKLT}{\mathop{\mathrm{MC}_\mathrm{non-klt}}}
\newcommand{\Diff}{\mathop{\mathrm{Diff}}}
\title[log quasi-numerically positive]{On log canonical divisors that are log quasi-numerically positive}
\author{Shigetaka Fukuda}
\address{Faculty of Education, Gifu Shotoku Gakuen University \\ Yanaizu-cho, Gifu 501-6194, Japan}
\email{fukuda@ha.shotoku.ac.jp}
\date{\empty}

\begin{document}
\begin{abstract}
Let $(X, \Delta)$ be a four-dimensional log variety that is projective over the field of complex numbers.
Assume that $(X, \Delta)$ is not Kawamata log terminal (klt) but divisorial log terminal (dlt).
First we introduce the notion of ``log quasi-numerically positive", by relaxing that of ``numerically positive".
Next we prove that, if the log canonical divisor $K_X + \Delta$ is log quasi-numerically positive on $(X, \Delta)$ then it is semi-ample.
\end{abstract}
\maketitle
\tableofcontents

\section{Introduction}

Throughout the paper every variety is projective over the field of complex numbers.
We follow the notation and terminology of the proceedings \cite{Utah} of ``the second Utah seminar''.

\begin{Defn}
A $\mathbf{Q}$-Cartier $\mathbf{Q}$-divisor $L$ on a projective variety $X$ is {\it numerically positive} ({\it nup}, for short) if $(L,C)>0$ for every curve $C$ on $X$.
A nef $\mathbf{Q}$-divisor $L$ on $X$ is {\it quasi-numerically positive} ({\it quasi-nup}, for short) if there exists a union $V$ of at most countably many Zariski-closed subsets $\subsetneqq X$ such that $(L,C)>0$ for every curve $C$ not contained in $V$.
A quasi-nup $\mathbf{Q}$-divisor $L$ on $X$ is {\it log quasi-numerically positive} ({\it log quasi-nup}, for short) on a divisorial log terminal (dlt) variety $(X, \Delta)$ if $L \vert_B$ is quasi-nup for every non-Kawamata log terminal (non-klt) center $B$ (in other words, for every $B \in \CNKLT (X, \Delta)$, under the notation of Section \ref{Sect:subad}).
\end{Defn}

Of course, the nupness (resp.~the log quasi-nupness) implies the log quasi-nupness (resp.~the quasi-nupness).
In the case where $(X, \Delta)$ is Kawamata log terminal (klt), the quasi-nupness is equivalent to the log quasi-nupness.

Recently F.~Ambro (\cite{Am}) reduced the famous log abundance conjecture, which claims that the nef log canonical divisors should be semi-ample, for klt varieties to the log minimal model conjecture and the problem of semi-ampleness of the quasi-nup log canonical divisors:

\begin{Prob}\label{Prob:A}
Assume that $(S,D)$ is klt and $K_S + D$ is quasi-nup.
Is $K_S + D$ semi-ample?
\end{Prob}

With regard to Problem \ref{Prob:A}, we consider the following

\begin{Prob}\label{Prob:B}
Assume that $(X,\Delta)$ is not klt but dlt and $K_X + \Delta$ is log quasi-nup on $(X,\Delta)$.
Is $K_X + \Delta$ semi-ample?
\end{Prob}

We note that the log abundance conjecture (including Problems \ref{Prob:A} and \ref{Prob:B}) for dlt varieties is known to be true in dimension $\leq 3$ (\cite{Ft}, \cite{KeMaMc}).
The subadjunction theory of Kawamata-Shokurov (in Section \ref{Sect:subad}) and a uniruledness theorem of Mori-Miyaoka type (due to Matsuki \cite{Ma}) enable us to reduce Problem \ref{Prob:B} (where $\lfloor \Delta \rfloor \neq 0$) in dimension $n$ to Problem \ref{Prob:A} (where $\lfloor D \rfloor = 0$) in dimension $\leq n-1$ (see Proposition \ref{Prop:Red}) and obtain the following main theorem.

\begin{Thm}\label{Thm:MT}
Assume that $(X,\Delta)$ is not klt but dlt, $\dim X =4$ and $K_X + \Delta$ is log quasi-nup on $(X,\Delta)$.
Then $K_X + \Delta$ is semi-ample.
\end{Thm}

In the case where $X$ is smooth and $\Delta$ is reduced and with only simple normal crossings, the theorem was proved in \cite{Fuk}.

\section{The subadjunction theory of Kawamata-Shokurov}\label{Sect:subad}
We recall the subadjunction theory of Kawamata (cf.~\cite[Lemma 5-1-9]{KMM}) and Shokurov (\cite[Subsection 3.2.3]{Sh}), clarify the notion of minimal non-klt centers and fix the relevant notation.

A log variety $(X,\Delta)$ consists of a normal variety $X$ and an effective $\mathbf{Q}$-divisor $\Delta$ on $X$ such that $\lceil \Delta \rceil$ is reduced.
The log variety $(X, \Delta)$ is said to be {\it divisorial log terminal} ({\it dlt}, for short) if $K_X + \Delta$ is $\mathbf{Q}$-Cartier and there exists a projective log resolution $f:Y \to X$ such that $K_Y + f^{-1}_* \Delta = f^* (K_X + \Delta) + E$ with the properties that $\lceil E \rceil \geq 0$, $\Exc (f)$ is a divisor and $\Supp (f^{-1}_* \Delta) \cup \Exc (f)$ is a simply normal crossing divisor.
(Moreover if $\lfloor \Delta \rfloor = 0$, the log variety $(X,\Delta)$ is said to be {\it Kawamata log terminal} ({\it klt}, for short).)

We set $D_i := f^{-1}_* \Delta_i$, where $\lfloor \Delta \rfloor = \sum_{i=1}^l \Delta_i$ and $\Delta_i$ is a prime divisor.
Define $\Strata_f (X,\Delta) := \{ \Gamma; \ k \geq 1, \ 1 \leq i_1 < i_2 < \ldots < i_k \leq l, \ \Gamma$ is an irreducible component of $D_{i_1} \cap D_{i_2} \cap \ldots \cap D_{i_k} \neq \emptyset \}$.
The set of non-klt centers $\CNKLT (X, \Delta) := \{ f(\Gamma) ; \ \Gamma \in \Strata_f (X,\Delta) \}$ is known not to depend on the choice of $f$.
Note that $\Exc (f) \nsupseteq \Gamma$ for every $\Gamma \in \Strata_f (X,\Delta)$, because $\Exc (f)$ is a divisor and $\Supp (\sum_{i=1}^l D_i) \cup \Exc (f)$ is a simply normal crossing divisor.
Thus for every $B \in \CNKLT (X, \Delta)$, the morphism $f$ is isomorphic over some suitable neighborhood of the generic point of $B$.
Therefore $f(\Gamma_1) \subset f(\Gamma_2)$ if and only if $\Gamma_1 \subset \Gamma_2$, for $\Gamma_1, \Gamma_2 \in \Strata_f (X,\Delta)$.
Consequently the set of minimal non-klt centers $\MCNKLT (X, \Delta) := \{ B ; \ B$ is a minimal element (with respect to inclusion) of $\CNKLT (X, \Delta) \}$ coincides with the set $\{ f(\Gamma) ; \ \Gamma$ is a minimal element of $\Strata_f (X,\Delta) \}$.

Now we focus on the subvariety $\Delta_i$ of $X$.
From Koll\'ar-Mori (\cite[Corollary 5.52]{KoMo}), $\Delta_i$ is normal and from the subadjunction theorem (Kawamata-Matsuda-Matsuki \cite[Lemma 5-1-9]{KMM}), $\Diff_{\Delta_i} (\Delta - \Delta_i) \geq 0$, where $(K_X + \Delta) \vert_{\Delta_i} = K_{\Delta_i} + \Diff_{\Delta_i} (\Delta - \Delta_i)$.
Put $f_i := f \vert_{D_i}$ and $\Strata_f (X,\Delta) \vert_{D_i} := \{ \Gamma \in \Strata_f (X,\Delta) ; \ \Gamma \subsetneqq D_i \}$.
Note that $K_{D_i} + f^{-1}_* (\Delta - \Delta_i) \vert_{D_i} = f_i^* (K_{\Delta_i} + \Diff_{\Delta_i} (\Delta - \Delta_i)) + E \vert_{D_i}$ and that $(f^{-1}_*(\Delta - \Delta_i) \vert_{D_i} - E \vert_{D_i}) = \sum \{ D_j \vert_{D_i} ; \ j \neq i, \ D_j \cap D_i \neq \emptyset \} + F$ (where $F$ is some divisor such that $\Supp F$ does not contain any irreducible component of $D_j \vert_{D_i} \neq 0$ and that $- \lfloor F \rfloor$ is effective and $f_i$-exceptional because $\Diff_{\Delta_i}(\Delta - \Delta_i) \geq 0$).
Here $\Exc (f_i) \nsupseteq \Gamma$ for any $\Gamma \in \Strata_f (X,\Delta) \vert_{D_i}$, since $\Exc (f_i) \subseteq \Exc (f) \cap D_i$.
Hence, by considering a suitable embedded resolution of $\Exc (f_i) \subseteq D_i$, we obtain that $(\Delta_i , \Diff_{\Delta_i} (\Delta - \Delta_i))$ is dlt (Shokurov \cite[Subsection 3.2.3]{Sh} and Fujino \cite[Proof of Theorem 0.1]{Fn}).

Then $\CNKLT (\Delta_i , \Diff_{\Delta_i} (\Delta - \Delta_i)) = \{f_i (\Gamma) ; \ \Gamma \in \Strata_f (X,\Delta) \vert_{D_i} \}$ and $\CNKLT(X,\Delta) = \bigcup_{j=1}^l (\CNKLT
 (\Delta_j, \Diff_{\Delta_j} (\Delta - \Delta_j)) \cup \{ \Delta_j \} )$.
Note that $f(\Gamma_1) \subset f(\Gamma_2) \subset \Delta_i$ if and only if $\Gamma_1 \subset \Gamma_2 \subset D_i$, for $\Gamma_1, \Gamma_2 \in \Strata_f (X,\Delta)$.
Therefore $\MCNKLT (\Delta_i , \Diff_{\Delta_i} (\Delta - \Delta_i)) = \MCNKLT (X, \Delta) \cap$ $\CNKLT (\Delta_i , \Diff_{\Delta_i} (\Delta - \Delta_i))$.

We define the maximal dimension of minimal non-klt centers by $l(X,\Delta) :=$ $ \max \{ \dim B ; \ B \in \MCNKLT (X,\Delta) \}$ in the case where $(X, \Delta)$ is not klt but dlt.
Of course, $l(X,\Delta) \leq \dim X - 1$ in this case.

\section{Reduction of the non-klt but dlt case, to the klt case in lower dimensions}

We reduce Problem \ref{Prob:B}, to Problem \ref{Prob:A} in lower dimensions.

\begin{Prop}\label{Prop:Red}
Let $(X,\Delta)$ be a log variety that is not klt but dlt and whose log canonical divisor $K_X + \Delta$ is log quasi-nup on $(X,\Delta)$.
Assume that Problem \ref{Prob:A} has an affirmative answer in dimension $\leq l(X,\Delta)$ {\rm (}the maximal dimension of minimal non-klt centers{\rm \/ )}.
Then $K_X + \Delta$ is semi-ample.
\end{Prop}

\begin{proof}
We shall prove the proposition by induction on $n:= \dim X$, heavily relying on the notation introduced in Section \ref{Sect:subad}.

Note that $(K_X + \Delta) \vert_{\Delta_i} = K_{\Delta_i} + \Diff_{\Delta_i} (\Delta - \Delta_i)$ is log quasi-nup on the dlt variety $(\Delta_i, \Diff_{\Delta_i} (\Delta - \Delta_i))$ from the fact that $(\CNKLT 
 (\Delta_i, \Diff_{\Delta_i} (\Delta - \Delta_i)) \cup \{ \Delta_i \}) \subseteq \CNKLT (X, \Delta)$.
Because $\MCNKLT (\Delta_i , \Diff_{\Delta_i} (\Delta- \Delta_i)) =$ $ \MCNKLT (X, \Delta) \cap \CNKLT (\Delta_i , \Diff_{\Delta_i} (\Delta- \Delta_i))$, the inequality $l(\Delta_i, \Diff_{\Delta_i} (\Delta - \Delta_i)) \leq l(X, \Delta)$ holds in the case where $(\Delta_i, \Diff_{\Delta_i} (\Delta - \Delta_i))$ is not klt. 
Therefore we know that $(K_X + \Delta) \vert_{\Delta_i}$ is semi-ample from the induction hypothesis in this case.
The $\mathbf{Q}$-divisor $(K_X + \Delta) \vert_{\Delta_i}$ is semi-ample also in the case where $(\Delta_i, \Diff_{\Delta_i} (\Delta - \Delta_i))$ is klt, because the value of $l(X,\Delta)$ becomes $n-1$ and hence the assumption of the theorem applies.
Thus $((K_X + \Delta) \vert_B)^{\dim B} > 0$ for every $B \in \CNKLT (X, \Delta) = \bigcup_{j=1}^l (\CNKLT
 (\Delta_j, \Diff_{\Delta_j} (\Delta - \Delta_j)) \cup \{ \Delta_j \} )$.

Next we show that $(K_X + \Delta)^n >0$.
By assuming that $(K_X + \Delta)^n = 0$, we shall imply the contradiction.
Note that $-K_Y f^* (K_X + \Delta)^{n-1} = (f^{-1}_* \Delta - E - f^* (K_X + \Delta)) f^* (K_X + \Delta)^{n-1} = (f^{-1}_* \Delta - E) f^* (K_X + \Delta)^{n-1} = (f^{-1}_* \Delta) f^* (K_X + \Delta)^{n-1} \geq (f^{-1}_* \Delta_1) f^* (K_X + \Delta)^{n-1} = ((K_X + \Delta) \vert_{\Delta_1})^{n-1} = ((K_X + \Delta) \vert_{\Delta_1})^{\dim \Delta_1} > 0$ because $E$ is $f$-exceptional and $\Delta_1 \in \CNKLT (X, \Delta)$.
Thus from Matsuki (the uniruledness theorem of Mori-Miyaoka type, \cite{Ma}), $Y$ is covered by $f^* (K_X + \Delta)$-trivial curves.
Therefore also $X$ is covered by $(K_X + \Delta)$-trivial curves.
This is a contradiction, because $K_X + \Delta$ is quasi-nup!
So we have that $(K_X + \Delta)^n >0$.

Consequently $K_X + \Delta$ becomes {\it nef and log big} on $(X, \Delta)$ (i.e.~$K_X + \Delta$ is nef, $(K_X + \Delta)^n >0$ and also $((K_X + \Delta) \vert_B )^{\dim B} >0$ for every $B \in \CNKLT (X, \Delta)$) and thus $K_X + \Delta$ is semi-ample, by virtue of the base point free theorem of Reid-type (Fujino \cite{Fn}).
\end{proof}

\section{Proof of the main theorem}

\begin{proof}[Proof of Theorem \ref{Thm:MT}]
Note that $l(X, \Delta) \leq \dim X -1 = 3$.
Thus the assumption of Proposition \ref{Prop:Red} is satisfied, from the log abundance theorem (\cite{Ft}, \cite{KeMaMc}) for klt varieties in dimensions $2$ and $3$.
\end{proof}


\begin{thebibliography}{99}

\bibitem{Am}
F.~Ambro, {\em The moduli b-divisor of an lc-trivial fibration},
math. AG/0308143, August 2003.

\bibitem{Fn}
O.~Fujino, {\em Base point free theorem of Reid-Fukuda type},
J. Math. Sci. Univ. Tokyo {\bf 7}(2000), 1--5, MR1749977.

\bibitem{Ft}
T.~Fujita, {\em Fractionally logarithmic canonical rings of algebraic surfaces},
J. Fac. Sci. Univ. Tokyo Sect. IA Math. {\bf 30}(1984), 685--696, MR0731524.

\bibitem{Fuk}
S.~Fukuda, {\em A note on the ampleness of numerically positive log canonical and anti-log canonical divisors},
math. AG/0305357, May 2003.

\bibitem{KMM}
Y.~Kawamata, K.~Matsuda and K.~Matsuki, {\em Introduction to the minimal model problem},
Algebraic geometry (Sendai, 1985), Adv. Stud. Pure Math., no. 10, North-Holland, Amsterdam, 1987, pp. 283--360, MR0946243.

\bibitem{KeMaMc}
S.~Keel, K.~Matsuki and J.~McKernan, {\em Log abundance theorem for threefolds},
Duke Math. J. {\bf 75}(1994), 99--119, MR1284817.

\bibitem{Utah}
J.~Koll\'ar et al., {\em Flips and abundance for algebraic threefolds},
Papers from the Second Summer Seminar on Algebraic Geometry held at the University of Utah (Salt Lake City, Utah, 1991), Ast\'erisque no. 211, Soci\'et\'e Math\'ematique de France, Paris, 1992, pp. 1--258, MR1225842.

\bibitem{KoMo}
J.~Koll\'ar and S.~Mori, {\em Birational geometry of algebraic varieties},
Cambridge Tracts in Mathematics, no. 134, Cambridge University Press, Cambridge, 1998, MR1658959.

\bibitem{Ma}
K.~Matsuki, {\em A correction to the paper ``Log abundance theorem for threefolds"},
math. AG/0302360, February 2003.

\bibitem{Sh}
V.~Shokurov, {\em 3-fold log flips},
Russian Acad. Sci. Izv. Math. {\bf 40}(1993), 95--202, MR1162635.

\end{thebibliography}
\end{document}